\documentclass[a4paper,12pt]{article}

\usepackage{amsfonts}
\usepackage{graphicx}
\usepackage{subfig}
\usepackage{amsmath}
\usepackage{amssymb}
\usepackage{enumitem}
\usepackage{url}
\usepackage{color}
\usepackage{authblk}
\usepackage{amsthm}

\newtheorem{proposition}{Proposition}
\newtheorem{theorem}{Theorem}
\newtheorem{corollary}{Corollary}
\newtheorem*{proof*}{Proof}
\newtheorem{remark}{Remark}

\title{\vspace{2cm}{\bf {\sc Similarity solutions for thawing processes with a convective boundary condition}}\vspace{1cm}}

\author{Ceretani, Andrea N. and Tarzia, Domingo A.\\{\small CONICET-Depto. Matem\'atica, FCE, Universidad Austral,
Paraguay 1950, S2000FZF Rosario, Argentina. E-mail: \textcolor{blue}{aceretani@austral.edu.ar}; \textcolor{blue}{dtarzia@austral.edu.ar}}}

\date{}

\begin{document}  
\maketitle

\thispagestyle{empty} 

\begin{abstract}
Similarity solutions for a one-dimensional mathematical model for thawing in a saturated semi-infinite porous media is considered when
change of phase induces a density jump and a convective boundary condition is imposed at the fixed face $x=0$. Different cases depending
on physical parameters are analysed and an explicit solution of a similarity type is obtained if and only if an inequality for data is
verified. Moreover, a monotony property respect to the coefficient which characterizes the heat transfer at the fixed face $x=0$ is
obtained for the coefficient involved in the definition of the free boundary.

Relationship between the Stefan problem with convective condition at $x=0$ considered in this paper and the Stefan problem with temperature
condition at the same face studied in (Fasano-Primicerio-Tarzia, Math. Models Meth. Appl. Sci., 9 (1999), 1-10) is analized and conditions for
physical parameters under which both problems became equivalents are obtained. Furthermore, an inequality to be satisfied for the coefficient
which characterizes the free boundary of the problem with a temperature or a convective boundary condition at the fixed face $x=0$ is also obtained.
\end{abstract}


{\bf 2010 AMS subjet classification}: 35R35 - 35C06 - 80A22

\newpage

\section{Introduction}

In this paper, we consider the problem of thawing of a semi-infinite partially frozen porous media saturated with an incompressible liquid
when change of phase induces a density jump and a convective boundary condition is imposed on the fixed face, with the
aim of constructing similarity solutions (for a detailed exposition of the physical background we refer to \cite{ChRu1992,FaGuPrRu1993,Na1990,OnMi1985,Ta1997}). In \cite{FaPrTa1999} and \cite{LoTa2001} (which generalized \cite{Ta1981-1982}) similarity solutions
are obtained when a temperature and a heat flux condition are imposed at the fixed boundary, respectively.
In this paper, we deal with the
same physical situations as in \cite{FaPrTa1999,LoTa2001} and we study a one-dimensional model of the problem where the unknowns are the
temperature $u(x,t)$ of the unfrozen zone $Q_1=\{(x,t): 0<x<s(t), t>0\}$, the temperature $v(x,t)$ of the frozen zone
$Q_2=\{(x,t): x>s(t), t>0\}$ and the free boundary $x=s(t)$, defined for $t>0$, separating $Q_1$ and $Q_2$, which satisfies the following
equations and boundary and initial conditions (we refer to \cite{FaGuPrRu1993} for a detailed explanation of the model):
\begin{align}
\label{1}&u_t=d_Uu_{xx}-b\rho\dot{s}(t)u_x	&0<x<s(t), \hspace{0.5cm}t>0 \\
\label{2}&v_t=d_Fv_{xx}	&x>s(t),\hspace{0.5cm}t>0\\
\label{3}&u(s(t),t)=v(s(t),t)=d\rho s(t)\dot{s}(t)	&	t>0\\
\label{4}&k_Fv_x(s(t),t)-k_Uu_x(s(t),t)=\alpha\dot{s}(t)+\beta\rho s(t)(\dot{s}(t))^2	&	t>0\\
\label{5}&v(x,0)=v(+\infty,t)=-A	&x>0,\hspace{0.5cm}t>0\\
\label{6}&s(0)=0\\
\label{7}&k_Uu_x(0,t)=\frac{h_0}{\sqrt{t}}(u(0,t)-B)	&  t>0
\end{align}
with:
\begin{equation*}
d_U=\alpha_U^2=\frac{k_U}{\rho_Uc_U}\text{,}
\hspace{1cm}
d_F=\alpha_F^2=\frac{k_F}{\rho_Fc_F}\text{,}
\hspace{1cm}
b=\frac{\epsilon\rho_Wc_W}{\rho_Uc_U}\text{,}
\hspace{1cm}
d=\frac{\epsilon\gamma\mu}{K}\text{,}
\end{equation*}
\begin{equation*}
\rho=\frac{\rho_W-\rho_I}{\rho_W}\text{,}
\hspace{1cm}
\alpha=\epsilon\rho_I l\text{,}
\hspace{1cm}
\beta=\frac{\epsilon^2\rho_I(c_W-c_I)\gamma\mu}{K}=\epsilon d\rho_I(c_W-c_I)\neq 0\text{.}
\end{equation*}
where:\\
$\epsilon$: porosity,\\
$\rho_W$ and $\rho_I$: density of water and ice ($g/cm^3$),\\
$c$: specific heat at constant density ($cal/g\,^\circ C$),\\
$k_F$ and $k_U$: conductivity of the frozen and unfrozen zones ($cal/s\,cm\,^\circ C$),\\
$u$: temperature of the unfrozen zone ($^\circ C$),\\
$v$: temperature of the frozen zone ($^\circ C$),\\
$u=v=0$: melting point at atmospheric pressure,\\
$l$: latent heat at $u=0$ ($cal/g$),\\
$\gamma$: coefficient in the Clausius-Clapeyron law ($s^2cm^\circ C/g$),\\
$\mu>0$: viscosity of the liquid ($g/s\,cm$),\\
$K>0$: hydraulic permeability ($cm^2$),\\
$B>0$: external boundary temperature at the fixed face $x=0$ ($^\circ C$),\\
$B_0>0$: temperature at the fixed face $x=0$ ($^\circ C$),\\
$-A<0$: initial temperature ($^\circ C$),\\
$h_0>0$: coefficient which characterizes the heat transfer at the fixed face $x=0$ ($Cal/s^\frac{1}{2}\,cm^2\,^\circ C$).

\begin{remark}
The free boundary problem (\ref{1})-(\ref{7}) reduces to the usual Stefan problem when
\begin{equation*}
\rho=0
\end{equation*}
since in that case we have the cassical Stefan conditions on $x=s(t)$, i.e.:
\begin{equation*}
u(s(t),t)=v(s(t),t)=0\,,\quad t>0
\end{equation*}
\begin{equation*}
k_Fv_x(s(t),t)-k_Uu_x(s(t),t)=\alpha\dot{s}(t)
\end{equation*}
and therefore from now on we assume that $\rho\neq 0$.
\end{remark}

The goal of this paper is to find the necessary and/or sufficient conditions for data (with three dimensionless parameters) in order to
obtain an instantaneous phase-change process (\ref{1})-(\ref{7}) with the corresponding explicit solution of the similarity type when a
convective boundary condition of type (\ref{7}) is imposed on the fixed face $x=0$ \cite{AlSo1993,CaJa1959,Ru1971,We1901,ZuCh1994}. 
We remark that the solution given in \cite{ZuCh1994} is not correct for any data (in particular, for small heat transfer coefficient) for the classical two-phase Stefan problem ($\rho=0$) which was improved in \cite{Ta2004} obtaining the necessary and sufficient condition to get the corresponding explicit solution. Furthermore, we study the relationship between the problem (\ref{1})-(\ref{7})
and the problem studied in \cite{FaPrTa1999}, which consists of equations (\ref{1})-(\ref{6}) and the following temperature boundary condition at
the fixed face $x=0$:
\begin{equation}
\label{8}u(0,t)=B_0\hspace{2cm}t>0\text{, }
\end{equation}
where $B_0$ is the boundary temperature at the fixed face $x=0$, with the aim of finding conditions under which both problems
become equivalents.

Recently Stefan-like problems were studied in \cite{BrNa2012,BrNa2014,ChKo2012,Ko2012,MiVy2014,RoSa2014,SaTa2011,VoFa2013,Vo2014,ZhWaBu2014}. The plan is the following: we first obtain in Section \ref{simSol} the necessary and sufficient condition in order to have a similarity
solution of the free boundary problem (\ref{1})-(\ref{7}) as a function of a positive parameter which must be the solution of a
transcendental equation with three dimensionless parameters defined by the thermal coefficients, and initial and boundary conditions. We also find a monotony property between the coefficients which characterize the free boundary and the
the heat transfer at the fixed face $x=0$.
In Section \ref{exisUniq} we give the necessary and/or sufficient conditions for the three real parameters involved in the trascendental
equation in order to obtain an
instantaneous phase-change process (\ref{1})-(\ref{7}) with the corresponding similarity solution. We generalize results obtained for particular cases given in \cite{NaSaTa2010,Ta2004}. Finally, in Section \ref{twoProb} we
analize the relationship between the problems (\ref{1})-(\ref{7}) and (\ref{1})-(\ref{6}) and (\ref{8}), and we obtain conditions for data under
which both problem became equivalents. Furthermore, we obtain an inequality which satisfies the coefficient involved in the definition
of the free boundary of the problem (\ref{1})-(\ref{6}) and (\ref{8}).

\section{Similarity solutions}\label{simSol}

We have:

\begin{theorem}\label{th:1}
The free boundary value problem (\ref{1})-(\ref{7}) has the similarity solution:
\begin{align}
\label{u}&u(x,t)=\frac{Bg(p,\xi)+\frac{AMk_U}{2h_0\alpha_U}\xi^2+(AM\xi^2-B)\displaystyle\int_0^{\frac{x}{2\alpha_U\sqrt{t}}}\text{ exp}(-r^2+pr\xi)\,dr}{g(p,\xi)+\frac{k_U}{2h_0\alpha_U}}\\
\label{v}&v(x,t)=\frac{AM\xi^2+A\text{ erf}(\gamma_0\xi)-A(1+M\xi^2)\text{ erf}\left(\frac{x}{2\alpha_F\sqrt{t}}\right)}{\text{erfc}{(\gamma_0\xi})}\\
\label{s}&s(t)=2\xi\alpha_U\sqrt{t}
\end{align}
if and only if the dimensionless coefficient $\xi>0$ satisfies the following equation:
\begin{equation}\label{eq:xi}
G(M,p,y)=y+Ny^3\text{,}\hspace{2cm}y>0
\end{equation}
involving three dimensionless parameters, $N$, $M$ and $p$, defined by:
\begin{equation}
N=\frac{2\beta\rho\alpha_U^2}{\alpha}\in\mathbb{R}\text{,}
\hspace{2cm}
M=\frac{2d\rho\alpha_U^2}{A}\in\mathbb{R}\text{,}
\hspace{2cm}
p=2b\rho\in\mathbb{R}\text{,}
\end{equation}
where:
\begin{flalign}
\label{G}&G(M,p,y)=\delta_1\left(1-\frac{AM}{B}y^2\right)G_1(p,y)-\delta_2(1+My^2)G_2(y)\hspace{1cm}y>0\text{,}\\
\label{G_1}&G_1(p,y)=\frac{\text{exp}((p-1)y^2)}{K_0+g(p,y)}\hspace{5.6cm}p\in\mathbb{R},\,y>0\text{,}\\
\label{G_2}&G_2(y)=\frac{\text{exp}(-\gamma_0^2y^2)}{\text{erfc}(\gamma_0y)}\hspace{6.3cm}M\in\mathbb{R},\,y>0\\
\label{g}&g(p,y)=\displaystyle\int_0^y\text{ exp}(-r^2+pry)\,dr\hspace{4.7cm}p\in\mathbb{R},\,y>0\text{,}\\
&\delta_1=\frac{k_UB}{2\alpha\alpha_U^2}>0\text{,}\hspace{0.25cm}\delta_2=\frac{k_FA}{\alpha\alpha_U\alpha_F\sqrt{\pi}}>0\text{,}
\hspace{0.25cm}K_0=\frac{k_U}{2\alpha_Uh_0}>0\text{,}\hspace{0.25cm}\gamma_0=\frac{\alpha_U}{\alpha_F}>0\text{.}
\end{flalign}
\end{theorem}

\begin{proof*}
We will follow the method introduced in \cite{FaPrTa1999,LoTa2001}. First of all, we note that the function $u$ defined by:
\begin{equation*}
u(x,t)=\Phi(\eta)\hspace{1cm}\text{ with }\eta=\frac{x}{2\alpha_U\sqrt{t}}
\end{equation*}
is a solution of equation (\ref{1}) if and only if $\Phi$ satisfies the following equation:
\begin{equation*}
\frac{1}{2}\Phi''(\eta)+\left(\eta-\frac{b\rho}{\alpha_U}\dot{s}(t)\sqrt{t}\right)\Phi(\eta)=0
\end{equation*}
Then, if we consider the function $x=s(t)$ defined as in (\ref{s}), for some $\xi>0$ to be determined, we obtain that:
\begin{equation}\label{Phi}
\Phi(\eta)=C_1+C_2\displaystyle\int_0^\eta\text{exp}(-r^2+2b\rho\xi r)\,dr
\end{equation}
where $\xi$, $C_1$ and $C_2$ are constant values.

Therefore, a solution of equation (\ref{1}) is given by:
\begin{equation}
u(x,t)=C_1+C_2\displaystyle\int_0^{x/2\alpha_U\sqrt{t}}\text{exp}(-r^2+2b\rho\xi r)\,dr\hspace{1cm}
\end{equation}
where $C_1$, $C_2$ and $\xi$ are constant values to be determined.

On the other hand, it is well known that:
\begin{equation}
v(x,t)=C_3+C_4\text{erf}\left(\frac{x}{2\alpha_F\sqrt{t}}\right)\text{,}
\end{equation}
is a solution of equation (\ref{2}), where $C_3$ and $C_4$ are constant values to be determined and
$\text{erf}$ is the error function defined by:
\begin{equation*}
\text{erf}(x)=\frac{2}{\sqrt{\pi}}\displaystyle\int_0^x\text{exp}(-r^2)\,dr
\end{equation*}

From conditions (\ref{3}), (\ref{5}) and (\ref{7}) we have that:
\begin{equation*}
C_1=\frac{Bg(2b\rho,\xi)+\frac{d\rho\alpha_Uk_U}{h_0}\xi^2}{g(2b\rho,\xi)+\frac{k_U}{2h_0\alpha_U}}\text{,}\hspace{2cm}
C_2=\frac{2d\rho\alpha_U^2\xi^2-B}{g(2b\rho,\xi)+\frac{k_U}{2h_0\alpha_U}}\text{,}
\end{equation*}
\begin{equation*}
C_3=\frac{2d\rho\alpha_U^2\xi^2+A\text{erf}\left(\gamma_0\xi\right)}{\text{erfc}\left(\gamma_0\xi\right)}\text{,}\hspace{2cm}
C_4=-\frac{2d\rho\alpha_U^2\xi^2+A}{\text{erfc}\left(\gamma_0\xi\right)}
\end{equation*}
and by imposing condition (\ref{4}), we obtain that the functions $s$, $u$ and $v$ defined above by (\ref{u})-(\ref{s}), are a solution of the free boundary value problem
(\ref{1})-(\ref{7}) if and only if the dimensionless parameter $\xi$ satisfies the equation (\ref{eq:xi}), and therefore the thesis holds.
\begin{flushright}
$\blacksquare$
\end{flushright}
\end{proof*}

\begin{remark}
 The similarity solution (\ref{Phi}) is a generalized Neumann solution for the classical Stefan problem (see \cite{AlSo1993}) when a
 convective boundary condition at the fixed face $x=0$ is imposed (see \cite{Ta2004}).
\end{remark}

The following result give us a relationship between the coefficient $\xi$ which characterizes the free boundary of the problem
(\ref{1})-(\ref{7}) and the coefficient $h_0>0$ which characterizes the heat transfer at the fixed face $x=0$.

\begin{corollary}\label{co:1}
If the free boundary value problem (\ref{1})-(\ref{7}) has a unique similarity solution of the type (\ref{u})-(\ref{s}) and the
dimensionless coefficient $\xi>0$ satisfies:
\begin{equation}
 0<\xi<\sqrt{\frac{B}{AM}}\text{,}
\end{equation}
then $\xi$ is a strictly increasing function of $h_0>0$.
\end{corollary}

\begin{proof*}
Let us assume that the free boundary value problem (\ref{1})-(\ref{7}) has the similarity solution (\ref{u})-(\ref{s}). Due to previous
theorem we know that the dimensionless coefficient $\xi$ satisfies equation (\ref{eq:xi}). Since the LHS of equation (\ref{eq:xi}) is a strictly increasing function of $h_0>0$, for all $y\in\left(
0,\sqrt{\frac{B}{AM}}\right)$, as well as its limit when $y$ tends to $0^+$:
$$\displaystyle\lim_{y\to 0^+}\left(\delta_1\left(1-\frac{AM}{B}y^2\right)G_1(p,y)-
\delta_2(1+My^2)G_2(y)\right)=\frac{\delta_1}{K_0}-\delta_2\text{,}$$
and the RHS of equation (\ref{eq:xi}) is a strictly increasing function of $y>0$, we can conclude that $\xi$ is a strictly
increasing function of $h_0>0$.
\begin{flushright}
$\blacksquare$
\end{flushright}
\end{proof*}

\section{Existence and uniqueness of similarity solutions}\label{exisUniq}

Due to previous theorem, in order to analize the existence and uniqueness of similarity solutions of the free boundary value problem
(\ref{1})-(\ref{7}), we can focus in the solvability of equation (\ref{eq:xi}). With this aim, we split our analysis into four cases
which correspond to four possible combinations of the signs of the dimensionless parameters $M$ and $N$.

To prove the following results we will use properties of functions $g$, $G_1$ and $G_2$ involved in equation (\ref{eq:xi}) which we
proved in the Appendix.

\begin{theorem}\label{th:2}
If $M>0$ and $N>0$ then:
\begin{enumerate}
\item If:
\begin{equation}\label{h_0-1}
h_0>\frac{A}{B}\frac{k_F}{\sqrt{\pi d_F}}
\end{equation}
then equation (\ref{eq:xi}) has at least one solution $\xi$, which satisfies $0<\xi<\sqrt{\frac{B}{AM}}$.\label{th:2-a}
\item\label{th:2-b} If $p\leq 1$ then equation (\ref{eq:xi}) has a solution $\xi$, which satisfies $0<\xi<\sqrt{\frac{B}{AM}}$, if and only if
(\ref{h_0-1}) holds. Moreover, when (\ref{h_0-1}) holds, equation (\ref{eq:xi}) has a unique solution $\xi$, which satisfies $0<\xi<\sqrt{\frac{B}{AM}}$.
\end{enumerate}
\end{theorem}

\begin{proof*}
Let $M>0$ and $N>0$.
\begin{enumerate}

\item We have proved in the Appendix that:
\begin{equation*}
\displaystyle\lim_{y\to+\infty} \left(1-\frac{AM}{B}y^2\right)G_1(p,y)=\left\{\begin{array}{ccc}
                                                                              0 & & p<2\\
                                                                              -\infty & & p\geq 2
                                                                               \end{array}\right.
\end{equation*}
and
\begin{equation*}
\displaystyle\lim_{y\to+\infty} (1+My^2)G_2(y)=+\infty\text{.}
\end{equation*}

Then:
\begin{equation*}
\displaystyle\lim_{y\to+\infty} \left[\delta_1\left(1-\frac{AM}{B}y^2\right)G_1(p,y)-
\delta_2(1+My^2)G_2(y)\right]=-\infty\text{.}
\end{equation*}
Furthermore:
\begin{equation*}
\displaystyle\lim_{y\to 0^+} \left[\delta_1\left(1-\frac{AM}{B}y^2\right)G_1(p,y)-
\delta_2(1+My^2)G_2(y)\right]=\frac{\delta_1}{K_0}-\delta_2\text{.}
\end{equation*}

Let us assume that {\footnotesize $h_0>\frac{A}{B}\frac{k_F}{\sqrt{\pi d_F}}$}. Since this inequality is equivalent to {\footnotesize $\frac{\delta_1}{K_0}-\delta_2>0$},
we have that the last limit is positive. Now taking into account that the RHS of (\ref{eq:xi}) is an increasing function from $0$ to $+\infty$, it follows that
equation (\ref{eq:xi}) has at least one positive solution.

Moreover, since:
{\begin{align*}
&\displaystyle\lim_{y\to 0^+} \delta_1\left(1-\frac{AM}{B}y^2\right)G_1(p,y)= \frac{\delta_1}{K_0}\,\text{,}\hspace{0.2cm}
\displaystyle\lim_{y\to +\infty} \delta_1\left(1-\frac{AM}{B}y^2\right)G_1(p,y)= -\infty\,\text{,}\\
&\displaystyle\lim_{y\to 0^+} \delta_2(1+My^2)G_2(y)= \delta_2\,\text{,}\hspace{1.7cm}
\displaystyle\lim_{y\to +\infty} \delta_2(1+My^2)G_2(y)= +\infty
\end{align*}}
and
$$\frac{\delta_1}{K_0}>\delta_2\text{,}$$
we have that the LHS of (\ref{eq:xi}) has positive zeros $q_1< q_2<\dotsm$. Thus we can find a
solution $\xi$ of (\ref{eq:xi}) such that $\xi<q_1$.
To prove that $0<\xi<\sqrt{\frac{B}{AM}}$, only remains to note that each $q_k$ satisfies $q_k<\sqrt{\frac{B}{AM}}$.

\item Assume that $p\leq 1$ and $h_0>\frac{A}{B}\frac{k_F}{\sqrt{\pi d_F}}$. We know that (see \cite{FaPrTa1999}):
\begin{equation*}
\frac{\partial}{\partial y}G_1(p,y)<0,\hspace{1cm} y>0\text{.}
\end{equation*}
Then:
{ \begin{align*}
\frac{\partial}{\partial y}\left(1-\frac{AM}{B}y^2\right)G_1(p,y)=
\frac{-2AMy}{B}G_1(p,y)+
\left(1-\frac{AM}{B}y^2\right)\frac{\partial}{\partial y}G_1(p,y)<0\text{,}\\
0<y<\sqrt{\frac{B}{AM}}\text{,}
\end{align*}}
which implies that $\left(1-\frac{AM}{B}y^2\right)G_1(p,y)$ is strictly decreasing in $\left(0,\sqrt{\frac{B}{AM}}\right)$. Furthermore,
$\left(1-\frac{AM}{B}y^2\right)G_1(p,y)$ is negative in $\left(\sqrt{\frac{B}{AM}},+\infty\right)$. On the other hand, we know from the
Appendix that $(1+My^2)G_2(y)$ is a positive strictly incresing function of $y>0$. Therefore, the LHS of
(\ref{eq:xi}) is strictly decreasing in $\left(0,\sqrt{\frac{B}{AM}}\right)$ and negative in $\left(\sqrt{\frac{B}{AM}},+\infty\right)$. We also have that the limit of the LHS of (\ref{eq:xi}) when $y$ tends to $0^+$ is positive
because we are assuming that $h_0>\frac{A}{B}\frac{k_F}{\sqrt{\pi\alpha_F}}$. Then, since RHS of (\ref{eq:xi}) is strictly increasing from $0$ to $+\infty$, we have that
equation (\ref{eq:xi}) has a unique solution $\xi\in\left(0,\sqrt{\frac{B}{AM}}\right)$.

Now, let us assume that (\ref{eq:xi}) has a solution $\xi\in\left(0,\sqrt{\frac{B}{AM}}\right)$. Suppose
$h_0$ does not verify the inequality (\ref{h_0-1}), that is $\frac{\delta_1}{K_0}-\delta_2\leq 0$. As LHS of (\ref{eq:xi}) is strictly
decreasing in $\left(0,\sqrt{\frac{B}{AM}}\right)$, follows that LHS of (\ref{eq:xi}) is negative in $\left(0,\sqrt{\frac{B}{AM}}\right)$.
Since RHS of (\ref{eq:xi}) is positive in $\left(0,\sqrt{\frac{B}{AM}}\right)$, we have a contradiction. Therefore
$h_0>\frac{A}{B}\frac{k_F}{\sqrt{\pi d_F}}$ and the thesis holds.
\end{enumerate}
\begin{flushright}
$\blacksquare$
\end{flushright}
\end{proof*}

\begin{remark}
The inequality (\ref{h_0-1}) was obtained in \cite{Ta2004} for the particular case $\rho=0$.
\end{remark}
The following corollary summarizes previous results.

\begin{corollary}\label{co:2}
If $M>0$, $N>0$, $p\leq 1$ and $h_0>\frac{A}{B}\frac{k_F}{\sqrt{\pi d_F}}$ then the free boundary value problem (\ref{1})-(\ref{7}) has
a unique similarity solution of the type (\ref{u})-(\ref{s}) and the dimensionless coefficient $\xi$ is a strictly increasing function
of the parameter $h_0$ on the interval $\left(\frac{A}{B}\frac{k_F}{\sqrt{\pi d_F}},+\infty\right)$.
\end{corollary}

\begin{theorem}\label{th:3}
If $M>0$ and $N<0$ then, if:
\begin{equation*}
h_0>\frac{A}{B}\frac{k_F}{\sqrt{\pi d_F}}\hspace{1cm}\text{and}\hspace{1cm}B<\frac{AM}{|N|}\text{,}
\end{equation*}
 then equation (\ref{eq:xi}) has at least one solution $\xi$, which satisfies $0<\xi<\sqrt{\frac{B}{AM}}$.
\end{theorem}

\begin{proof*} Let {\small $M>0$} and {\small $N<0$}. We have that LHS of (\ref{eq:xi}) is positive in $(0,q_1)$,
where {\small $q_1<\sqrt{\frac{B}{AM}}$}. On the other hand, RHS of (\ref{eq:xi}) is positive in $\left(0,\sqrt{\frac{1}{|N|}}\right)$. Since
$\sqrt{\frac{B}{AM}}<\sqrt{\frac{1}{|N|}}$, we can find a solution $\xi$ of (\ref{eq:xi}) which satisfies $0<\xi<\sqrt{\frac{B}{AM}}$.
\begin{flushright}
$\blacksquare$
\end{flushright}
\end{proof*}

\begin{theorem}\label{th:4}
If $M<0$ and $N>0$ then:
\begin{enumerate}
\item If:
\begin{equation*}
h_0>\frac{A}{B}\frac{k_F}{\sqrt{\pi d_F}}
\end{equation*}
and the LHS of (\ref{eq:xi}) has positive zeros being $q_1$ the smaller one,
then equation (\ref{eq:xi}) has at least one solution $\xi$, which satisfies $0<\xi<q_1$.
\item If:
\begin{equation*}
h_0\leq \frac{A}{B}\frac{k_F}{\sqrt{\pi d_F}}\hspace{1cm}\text{and}\hspace{1cm}N<\delta_2\sqrt{\pi}|M|\gamma_0\text{,}
\end{equation*}
then equation (\ref{eq:xi}) has at least one positive solution.
\end{enumerate}
\end{theorem}

\begin{proof*}
Let $M<0$ and $N>0$.
\begin{enumerate}
 \item It is analogous to the proof of the Theorem \ref{th:2} (\ref{th:2-a}).
 \item Let us assume that $h_0\leq \frac{A}{B}\frac{k_F}{\sqrt{\pi d_F}}$ and $N<\delta_2\sqrt{\pi}|M|\gamma_0$.

 From the assymptotic behavior of $\left(1-\frac{AM}{B}y^2\right)G_1(p,y)$ and $(1+My^2)G_2(y)$ when $y$ tends to $+\infty$ (see
 Appendix), we have that:
{\small \begin{equation*}
 \displaystyle\lim_{y\to +\infty}\frac{\delta_1\left(1-\frac{AM}{B}y^2\right)G_1(p,y)-\delta_2(1+My^2)G_2(y)}{y+Ny^3}=
 \left\{\begin{array}{ccc}
 -\frac{\delta_2\sqrt{\pi}\gamma_0M}{N}& &p\leq 2\\
 -\frac{AM(p-2)}{BN}-\frac{\delta_2\sqrt{\pi}\gamma_0M}{N}& &p>2
 \end{array}\right.\text{.}
 \end{equation*}}
 Therefore:
 \begin{equation*}\label{li:7}
 \displaystyle\lim_{y\to +\infty}\frac{\delta_1\left(1-\frac{AM}{B}y^2\right)G_1(p,y)-\delta_2(1+My^2)G_2(y)}{y+Ny^3}
 >\frac{\delta_2\sqrt{\pi}\gamma_0|M|}{N}\text{.}
 \end{equation*}
 Then, it follows from the last inequality that the LHS of (\ref{eq:xi}) tends to $+\infty$ faster than the RHS when $y$ tends
 to $+\infty$. Now, taking into account that:
 \begin{equation*}
 \displaystyle\lim_{y\to 0^+} \delta_1\left(1-\frac{AM}{B}y^2\right)G_1(p,y)-\delta_2(1+My^2)G_2(y)=\frac{\delta_1}{K_0}-\delta_2<0\text{,}
 \end{equation*}
 we have that equation (\ref{eq:xi}) has at least one positive solution.
\end{enumerate}
\begin{flushright}
$\blacksquare$
\end{flushright}
\end{proof*}

\begin{remark}
When we have the condition $h_0\leq\frac{A}{B}\frac{k_F}{\sqrt{\pi d_F}}$ in the classical Stefan problem with $\rho=0$ we only obtain a heat transfer problem without a phase-change process \cite{Ta2004}.
\end{remark}

\begin{theorem}\label{th:5}
If $M<0$ and $N<0$ then:
\begin{enumerate}
\item If:
\begin{equation*}
h_0>\frac{A}{B}\frac{k_F}{\sqrt{\pi d_F}}
\end{equation*}
and the LHS of the equation (\ref{eq:xi}) has positive zeros being $q_1$ the smaller one which satisfies:
\begin{enumerate}
\item $q_1<\sqrt{1/|N|}$, then equation (\ref{eq:xi}) has at least two positive solutions, one of them satisfies $0<\xi<\sqrt{1/|N|})$.
\item $q_1=\sqrt{1/|N|}$, then equation (\ref{eq:xi}) has $\xi=q_1$ as solution.
\end{enumerate}
\item If:
\begin{equation*}
h_0<\frac{A}{B}\frac{k_F}{\sqrt{\pi d_F}}
\end{equation*}
then equation (\ref{eq:xi}) has at least one positive solution $\xi$.
\end{enumerate}
\end{theorem}

\begin{proof*}
Let $M<0$ and $N<0$.
\begin{enumerate}
 \item
 \begin{enumerate}
   \item Assume that $h_0>\frac{A}{B}\frac{k_F}{\sqrt{\pi d_F}}$.

  We have proved in the Appendix that:
  \begin{equation*}
   \displaystyle\lim_{y\to+\infty} \left(1-\frac{AM}{B}y^2\right)G_1(p,y)=\left\{\begin{array}{ccc}
                                                                              0 & & p<2\\
                                                                              +\infty & & p\geq 2
                                                                               \end{array}\right.
  \end{equation*}																																							
  and
	\begin{equation*}
	\displaystyle\lim_{y\to+\infty} (1+My^2)G_2(y)=-\infty\text{.}
  \end{equation*}
  Then:
   \begin{equation*}
     \displaystyle\lim_{y\to +\infty} \left[\delta_1\left(1-\frac{AM}{B}y^2\right)G_1(p,y)-\delta_2(1+My^2)G_2(y)\right]=+\infty\text{.}
   \end{equation*}
  Furthermore:
   \begin{equation*}
     \displaystyle\lim_{y\to 0^+} \left[\delta_1\left(1-\frac{AM}{B}y^2\right)G_1(p,y)-\delta_2(1+My^2)G_2(y)\right]=\frac{\delta_1}{K_0}-\delta_2\text{.}
   \end{equation*}
  On the other hand, the RHS of (\ref{eq:xi}) is positive in $\left(0,\sqrt{\frac{1}{|N|}}\right)$ and negative in
   $\left(\sqrt{\frac{1}{|N|}},+\infty\right)$. Then, since $\frac{\delta_1}{K_0}-\delta_2>0$ and $q_1<\sqrt{1/|N|}$ we have that equation
  (\ref{eq:xi}) has at least two positive solutions, one of them satisfies $0<\xi<\sqrt{\frac{1}{|N|}}$.
 \item It is inmediate.
\end{enumerate}
 \item Now, let us assume that $h_0<\frac{A}{B}\frac{k_F}{\sqrt{\pi d_F}}$. We have:
\begin{equation*}
\displaystyle\lim_{y\to +\infty} \left[\delta_1\left(1-\frac{AM}{B}y^2\right)G_1(p,y)-\delta_2(1+My^2)G_2(y)\right]=+\infty\,\text{.}
\end{equation*}
Furthermore:
\begin{equation*}
\displaystyle\lim_{y\to 0^+} \left[\delta_1\left(1-\frac{AM}{B}y^2\right)G_1(p,y)-\delta_2(1+My^2)G_2(y)\right]=
\frac{\delta_1}{K_0}-\delta_2\text{.}
\end{equation*}
Taking into account that $\displaystyle\lim_{y\to 0^+} \left(y+Ny^3\right)=0$ and
$\displaystyle\lim_{y\to +\infty} \left(y+Ny^3\right)=-\infty$, since {\footnotesize $\frac{\delta_1}{K_0}-\delta_2<0$} we have that equation
(\ref{eq:xi}) has at least one positive solution $\xi$.
\end{enumerate}
\begin{flushright}
$\blacksquare$
\end{flushright}
\end{proof*}

\begin{remark}
Previous results imply that under the conditions specified in each case it is posible to find a similarity solution of the problem
(\ref{1})-(\ref{7}). Moreover, that solution guarantees a change of phase, that is, $AM\xi^2<u(0,t)<B$.
\end{remark}

\section{Relationship between the solutions of the Stefan problem with convective and temperature boundary conditions}\label{twoProb}

In this section, we analyze the relationship between problem (\ref{1})-(\ref{7}) and problem (\ref{1})-(\ref{6}),(\ref{8}) studied in \cite{FaPrTa1999}, corresponding to a temperature condition at the fixed face $x=0$.

In \cite{FaPrTa1999} it was proved that if $M>0$, $N>0$ and $p\leq2$, then the problem (\ref{1})-(\ref{6}),(\ref{8}) has a unique
similarity solution of the type:
\begin{align}
\label{U}&U(x,t)=B_0+\frac{AM\omega^2-B_0}{g(p,y)}\displaystyle\int_0^{\frac{x}{2\alpha_U\sqrt{t}}}\text{ exp}(-r^2+pr\omega)\,dr\\
\label{V}&V(x,t)=\frac{AM\omega^2\text{erf}\left(\frac{x}{2\alpha_F\sqrt{t}}\right)+A\left(\text{ erf}(\gamma_0\omega)-
\text{erf}\left(\frac{x}{2\alpha_F\sqrt{t}}\right)\right)}{\text{erfc}{(\gamma_0\omega})}\\
\label{S}&S(t)=2\omega\alpha_U\sqrt{t}
\end{align}

where $\omega$ is the unique positive solution of the trascendental equation:
\begin{equation}\label{eq:omega}
\widetilde{G}(M,p,y)=y+Ny^3\text{,}
\end{equation}
with:
{\begin{align}
\label{Gtilde}&\widetilde{G}(M,p,y)=\tilde{\delta}_1\left(1-\frac{AM}{B_0}y^2\right)\widetilde{G_1}(p,y)-\delta_2(1+My^2)G_2(M,y)\hspace{0.5cm}y>0\text{,}\\
&\widetilde{G}_1(p,y)=\frac{\exp{(p-1)y^2}}{g(p,y)}\hspace{5.8cm} p\in\mathbb{R},\,y>0\text{,}\\
&\tilde{\delta}_1=\frac{k_UB_0}{2\alpha\alpha_1^2}\text{.}
\end{align}}

Furthermore, $0<\omega<\sqrt{\frac{B_0}{AM}}$.

Henceforth, we will only deal with situations in which existence and uniqueness of similarity solutions of type (\ref{u})-(\ref{s}) for
problem (\ref{1})-(\ref{8}) or of type (\ref{U})-(\ref{S}) for problem (\ref{1})-(\ref{6}),(\ref{8}), are guarantee.

\begin{theorem}\label{th:6}
If $M>0$, $N>0$, $p\leq 1$ and $h_0>\frac{A}{B}\frac{k_F}{\sqrt{\pi d_F}}$ then the dimensionless coefficient $\xi$ which characterizes
the free boundary of the problem (\ref{1})-(\ref{7}) satisfies:
\begin{equation}
0<\xi(h_0)<\omega_\infty\hspace{2cm}\forall\,h_0\in \left(\frac{A}{B}\frac{k_F}{\sqrt{\pi d_F}},+\infty\right)
\end{equation}
where $\omega_\infty$ is the coefficient which characterizes the free boundary of the problem (\ref{1})-(\ref{6}),(\ref{8}) when the
temperature condition is given by $B$.
\end{theorem}

\begin{proof*}
Assume that $M>0$, $N>0$, $p\leq 1$ and $h_0>\frac{A}{B}\frac{k_F}{\sqrt{\pi d_F}}$. We know from Corollary \ref{co:2} that the
dimensionless coefficient $\xi$ is a strictly increasing function of the coefficient $h_0$ on the interval
$\left(\frac{A}{B}\frac{k_F}{\sqrt{\pi d_F}},+\infty\right)$. We also know that $\xi$ satisfies equation (\ref{eq:xi}) which became:
\begin{equation}\label{eq:omega_infty}
\delta_1\left(1-\frac{AM}{B_0}y^2\right)\widetilde{G_1}(p,y)-\delta_2(1+My^2)G_2(M,y)=y+Ny^3
\end{equation}
when $h_0$ tends to $+\infty$. Only remains to note that this last equation (\ref{eq:omega_infty}) has a unique solution $\omega_\infty$ because (\ref{eq:omega_infty}) is the corresponding equation to problem
(\ref{1})-(\ref{6}),(\ref{8}) when $B_0=B$, which has a unique similarity solution under the hypothesis considered here.\hspace{5cm}$\blacksquare$
\end{proof*}

\begin{theorem}\label{th:7}
If $M>0$, $N>0$, $p\leq1$ and $h_0>\frac{A}{B}\frac{k_F}{\sqrt{\pi d_F}}$ then the similarity solution (\ref{u})-(\ref{s}) of the
the problem (\ref{1})-(\ref{7}), coincide with the similarity solution (\ref{U})-(\ref{S}) of the problem (\ref{1})-(\ref{6}),(\ref{8}),
when the external boundary temperature at $x=0$ is defined as:
\begin{equation}\label{B_0}
B_0=\frac{2h_0\alpha_UBg(p,\xi)+AMk_U\xi^2}{2h_0\alpha_Ug(p,\xi)+k_U}
\end{equation}
Moreover, the dimensionless parameters $\xi$ and $\omega$ which characterize the free boundary in each problem, are equals.
\end{theorem}

\begin{remark}
Note that $B_0$ is positive since $B_0=\frac{2h_0\alpha_UBg(p,\xi)+AMk_U\xi^2}{2h_0\alpha_Ug(p,\xi)+k_U}=u(0,t)$, $t>0$, and $u(0,t)>0$.
\end{remark}

\begin{proof*}
Assume that $M>0$, $N>0$, $p\leq1$, $h_0>\frac{A}{B}\frac{k_F}{\sqrt{\pi d_F}}$ and $B_0$ is defined as in (\ref{B_0}).
First of all, we note that function $\widetilde{G}(M,p,y)$, given in (\ref{Gtilde}), can be written as:
\begin{equation*}
\widetilde{G}(M,p,y)=\delta_1\left(\frac{B_0}{B}-\frac{AM}{B}y^2\right)\widetilde{G_1}(p,y)-\delta_2(1+My^2)G_2(y)\hspace{2cm}y>0\text{.}
\end{equation*}
Because of the definition of $B_0$, given in (\ref{B_0}), we have:
\begin{equation*}
\frac{B_0}{B}-\frac{AM}{B}y^2=\frac{g(p,\xi)\left(1-\frac{AM}{B}y^2\right)+\frac{AM}{B}K_0(\xi^2-y^2)}{g(p,\xi)+K_0}\text{.}
\end{equation*}
Then, equation (\ref{eq:omega}) can be written as:
{\small \begin{equation}\label{eq:omega.1}
\delta_1\left(\frac{g(p,\xi)}{g(p,y)}\left(1-\frac{AM}{B}y^2\right)+
\frac{AMK_0}{Bg(p,y)}(\xi^2-y^2)\right)G_1(p,y)-
\delta_2(1+My^2)G_2(y)=y+Ny^3\text{.}
\end{equation}}
Now, since $\xi$ satisfies equation (\ref{eq:xi}), it follows that $\xi$ satisfies equation (\ref{eq:omega.1}). Therefore, $\xi$
coincide with the coefficient $\omega$ which characterize the free boundary of the problem (\ref{1})-(\ref{8}). Finally, from elementary
calculations, we have that similarity solutions given in (\ref{u})-(\ref{s}) and (\ref{U})-(\ref{S}) are coincident.
\begin{flushright}
$\blacksquare$
\end{flushright}
\end{proof*}

In the same way that the proof of the previous theorem, it can be shown the following result.

\begin{theorem}\label{th:8}
If $M>0$, $N>0$ and $p\leq1$ then the similarity solution (\ref{U})-(\ref{S}) of the problem (\ref{1})-(\ref{6}),(\ref{8}) coincide with
the similarity solution (\ref{u})-(\ref{s}) of the problem (\ref{1})-(\ref{7}), when the coefficient which characterizes the heat transfer at the fixed face $x=0$ is defined by:
\begin{equation}\label{h_0-2}
h_0=\frac{k_U(B_0-AM\omega^2)}{2\alpha_U(B-B_0)g(p,\omega)}
\end{equation}
and the boundary temperature $B$ at $x=0$ is such that $B>B_0$.

Moreover, the dimensionless parameters $\omega$ and $\xi$ which define the free boundary in each problem, are equals.
\end{theorem}

\begin{remark}
Note that $h_0$ is a positive number since {\footnotesize $B_0-AM\omega^2>0$} and {\footnotesize $B-B_0>0$}.
\end{remark}

We can conclude now that in the sense established in Theorems \ref{th:7} and \ref{th:8}, Stefan problems with convective and temperature
conditions at the fixed face $x=0$ given by (\ref{1})-(\ref{7}) and (\ref{1})-(\ref{6}),(\ref{8}), respectively, are equivalent when inequality (\ref{h_0-1}) is verified by data.

\begin{theorem}\label{th:9}
If $M>0$, $N>0$ and $p\leq1$ then the coefficient $\omega$ which characterizes the free boundary of the problem (\ref{1})-(\ref{6}),(\ref{8}) satisfies the following inequality:
\begin{equation}\label{inq:omega}
\frac{B_0-AM\omega^2}{g(p,\omega)}>\frac{2\alpha_U k_F A(B-B_0)}{\alpha_F k_U B\sqrt{\pi}}\text{,}\hspace{1cm}\forall\,B>B_0
\end{equation}
\end{theorem}

\begin{proof*}
Assume that $M>0$, $N>0$ and $p\leq1$, and let $B>B_0$. We know from Theorem \ref{th:8}, that the problem (\ref{1})-(\ref{7}) has a similarity solution of the type (\ref{u})-(\ref{s}) when the external boundary temperature is given by $B$ and the coefficient $h_0$ is defined as in (\ref{h_0-2}). We also know that the coefficients which characterize the free boundary in problems (\ref{1})-(\ref{7}) and (\ref{1})-(\ref{6}), (\ref{8}) are equals, that is $\xi=\omega$. Therefore, since $0<\omega<\sqrt{\frac{B_0}{AM}}$ and $B>B_0$, we have that $0<\xi<\sqrt{\frac{B}{AM}}$. Then, due to Theorem \ref{th:2-b} inequality (\ref{h_0-1}) holds. Only remains to note that inequality (\ref{h_0-1}) becames inequality (\ref{inq:omega}) when $h_0$ is defined as in (\ref{h_0-2}).
\begin{flushright}
$\blacksquare$
\end{flushright}
\end{proof*}

By taking limit when $B$ tends to $+\infty$ into both sides of inequality (\ref{inq:omega}), we have the following corollary.

\begin{corollary}
If $M>0$, $N>0$ and $p\leq1$ then the coefficient $\omega$ which characterizes the free boundary of the problem (\ref{1})-(\ref{6}),(\ref{8}) satisfies the following inequality:
\begin{equation}\label{ineq:omega-1}
\frac{B_0-AM\omega^2}{g(p,\omega)}>\frac{2\alpha_UAk_F}{\alpha_Fk_U\sqrt{\pi}}
\end{equation}
\end{corollary}

\begin{remark}
For the classical Stefan problem with $\rho=0$ the inequality (\ref{ineq:omega-1}) for the coefficient $\omega$, which characterizes the free boundary $s(t)$, given by (\ref{S}), becames as:
\begin{equation*}
\text{ erf}({\omega})<\frac{B_0}{A}\frac{k_U}{k_F}\sqrt{\frac{d_F}{d_U}}
\end{equation*}
which was obtained in \cite{Ta1981-1982}.
\end{remark}

\section*{Acknowledgment}
This paper has been partially sponsored by the Projects PIP No. 0534 from CONICET-Univ. Austral (Rosario, Argentina).

\section{Appendix}

\begin{proposition}\label{pr:1}
For any $M\in\mathbb{R}$ and $p\in\mathbb{R}$:
\begin{enumerate}
\item If $p<2$ then:
\begin{equation}\label{li:1}
\displaystyle\lim_{y\to+\infty} \left(1-\frac{AM}{B}y^2\right)G_1(p,y)=0\text{.}
\end{equation}
\item If $p\geq 2$ then:
\begin{equation}\label{li:2}
\displaystyle\lim_{y\to+\infty} \left(1-\frac{AM}{B}y^2\right)G_1(p,y)=\left\{\begin{array}{ccc}
                                                                              -\infty & & M>0\\
                                                                              +\infty & & M<0
                                                                               \end{array}\right.
\end{equation}
and
\begin{equation}\label{li:3}
\left(1-\frac{AM}{B}y^2\right)G_1(p,y)\simeq \left\{\begin{array}{ccc}
			                            -\frac{2AM}{B\sqrt{\pi}}y^2 & &p=2\\
			                            -\frac{AM(p-2)}{B}y^3 & &p>2
                                                    \end{array}\right.\text{, when $y\to +\infty$.}
\end{equation}
\end{enumerate}
\end{proposition}

\begin{proof*}
Let $M$ and $p$ real numbers.

\begin{enumerate}
\item Proof of (\ref{li:1}) is an immediate consecuence of the following facts:
\begin{equation}\label{li:4-5}
\displaystyle\lim_{y\to+\infty}G_1(p,y)=0\hspace{1cm}\text{ and }\hspace{1cm} \displaystyle\lim_{y\to+\infty}y^2G_1(p,y)=0
\end{equation}

Then, we will prove (\ref{li:4-5}). We split the proof into three cases depending on the sign of the parameter $p$.
\begin{enumerate}
\item $0<p\leq1$\\
Since $\displaystyle\lim_{y\to+\infty}g(p,y)=+\infty$ (see \cite{FaPrTa1999}), we have:
\begin{equation*}
\displaystyle\lim_{y\to+\infty}G_1(p,y)=0
\end{equation*}

On the other hand,
{\small \begin{equation*}
\displaystyle\lim_{y\to+\infty}\left(y^2G_1(p,y)\right)^{-1}=
\displaystyle\lim_{y\to+\infty}\left(\frac{K_0}{y^2\text{exp}((p-1)y^2)}+\frac{g(p,y)}{y^2\text{exp}((p-1)y^2)}\right)=+\infty
\end{equation*}}
since $\displaystyle\lim_{y\to+\infty}\frac{g(p,y)}{y^2}=+\infty$ (see \cite{LoTa2001}). Then, $\displaystyle\lim_{y\to+\infty}y^2G_1(p,y)=0$.

\item $1<p<2$

First of all, we note that:
\begin{equation*}\label{g.1}
g(p,y)=\frac{\sqrt{\pi}}{2}\text{exp}\left(\left(\frac{p}{2}y^2\right)^2\right)\left(\text{erf}\left(\frac{p}{2}y\right)+
\text{erf}\left(\frac{2-p}{2}y\right)\right)\text{.}
\end{equation*}
Then:
{\small \begin{equation*}
\displaystyle\lim_{y\to+\infty}G_1(p,y)=\displaystyle\lim_{y\to+\infty}\frac{\text{exp}\left(-\left(\frac{p}{2}-1\right)^2y^2\right)}
{K_0\text{exp}\left(-\left(\frac{p}{2}y\right)^2\right)+\frac{\sqrt{\pi}}{2}}\left(\text{erf}\left(\frac{p}{2}y\right)+
\text{erf}\left(\frac{2-p}{2}y\right)\right)=0\text{.}
\end{equation*}}
The proof of $\displaystyle\lim_{y\to+\infty}y^2G_1(p,y)=0$ is similar to the previous one.

\item $p\leq0$

It is useful to write $G_1(p,y)$ in the following way:
{\small \begin{equation*}
G_1(p,y)=\left[K_0\text{exp}((1-p)y^2)+\frac{\sqrt{\pi}}{2}\text{exp}\left(\left(1-\frac{p}{2}\right)^2y^2\right)
\left(\text{erf}\left(\left(1-\frac{p}{2}\right)y\right)-\text{erf}\left(-\frac{p}{2}y\right)\right)\right]^{-1}\text{.}
\end{equation*}}
Then, since:
{\small \begin{equation*}
\displaystyle\lim_{y\to+\infty}\left[K_0\text{exp}((1-p)y^2)+\frac{\sqrt{\pi}}{2}\text{exp}\left(\left(1-\frac{p}{2}\right)^2y^2\right)
\left(\text{erf}\left(\left(1-\frac{p}{2}\right)y\right)-\text{erf}\left(-\frac{p}{2}y\right)\right)\right]=+\infty\text{,}
\end{equation*}}
we have that $\displaystyle\lim_{y\to+\infty}G_1(p,y)=0$.

The proof of $\displaystyle\lim_{y\to+\infty}y^2G_1(p,y)=0$ is similar to the previous one.
\end{enumerate}

\item Now, we have:
\begin{equation}\label{li:6}
\displaystyle\lim_{y\to+\infty} y^2G_1(p,y)=0\text{.}
\end{equation}
In fact:
{\small\begin{equation*}
\displaystyle\lim_{y\to+\infty} G_1(p,y)=
\displaystyle\lim_{y\to+\infty} \frac{2(p-1)y\text{exp}((p-1)y^2)}{\frac{\partial}{\partial y}g(p,y)}\geq
\displaystyle\lim_{y\to+\infty} \frac{2(p-1)y\text{exp}((p-1)y^2)}{\text{exp}((p-1)y^2)+pyg(p,y)}
\end{equation*}}
and
{\small\begin{equation*}
\displaystyle\lim_{y\to+\infty} \left(\frac{2(p-1)y\text{exp}((p-1)y^2)}{\text{exp}((p-1)y^2)+pyg(p,y)}\right)^{-1}=
\displaystyle\lim_{y\to+\infty} \frac{1}{2(p-1)y}+\frac{1}{2y(p-1)}\frac{yg(p,y)}{\text{exp}((p-1)y^2)}=0\text{.}
\end{equation*}}

For the last limit, we are using the fact that $\displaystyle\lim_{y\to+\infty} \frac{yg(p,y)}{\text{exp}((p-1)y^2)}=\frac{1}{p-2}$
(see \cite{LoTa2001}). Then $\displaystyle\lim_{y\to+\infty} y^2G_1(p,y)=0$.

Finally, (\ref{li:3}) follows from similar arguments.
\end{enumerate}
\begin{flushright}
$\blacksquare$
\end{flushright}
\end{proof*}

The following result is proved in \cite{FaPrTa1999}.

\begin{proposition}\label{pr:2}
We have:
\begin{equation}
\displaystyle\lim_{y\to+\infty} (1+My^2)G_2(y)=\left\{\begin{array}{ccc}
                                                        +\infty & & M>0\\
                                                        -\infty & & M<0
                                                        \end{array}\right.
\end{equation}
and
\begin{equation}
(1+My^2)G_2(y)\simeq \sqrt{\pi}\gamma_0My^3\text{ as $y\to+\infty$}.
\end{equation}

Furthermore, when $M>0$, $(1+My^2)G_2(y)$ is an increasing function of $y>0$.
\end{proposition}

\end{document}